\documentclass{tran-l}

 \theoremstyle{definition}
 
 \theoremstyle{remark}
 
 \numberwithin{equation}{subsection}

 \newcommand{\To}{\longrightarrow}

\begin{document}

\title[$P_{Z}(S)$-Metrics and $P_{Z}(S)$-Metric Spaces]
 {$P_{Z}(S)$-Metrics and $P_{Z}(S)$-Metric Spaces}

\author{Reza Jafarpour-Golzari}

\address{Department of Mathematics, Payame Noor University
,P.O.BOX 19395-3697 Tehran, Iran}
\address{Faculty of Mathematics, Institute for Advanced Studies
in Basic Science (IASBS), P.O.Box 45195-1159, Zanjan, Iran}

\email{r.golzary@iasbs.ac.ir}

\thanks{}



\date{}

\dedicatory{}

\commby{}


\begin{abstract}
In this paper, we define notions of
$P_{Z}(S)$-metric and $P_{Z}(S)$-metric space and we show that
every $P_{Z}(S)$-metric Space, analogous to an ordinary metric space
and generally, a $\Lambda$-metric space, is a topological space,
and in continuance, we show that from a topological point of view, some
properties of $P_{Z}(S)$-metric spaces, and $\Lambda$-metric
spaces, have coordination.
\end{abstract}

\maketitle

\section{Introduction}

A metric on a nonempty set $X$, is defined as $d\colon X\times X\To
\mathbb{R}$ which satisfies three following condition:

1) $\forall x,y\in X,\ \ 0\leq d(x,y),\
d(x,y)=0\Longleftrightarrow x=y$,

2) $\forall x,y\in X,\ d(x,y)=d(y,x)$,

3) $\forall x,y,z\in X,\ d(x,y)\leq d(x,z)+d(z,y)$.

The second condition, is known as symmetry and the third, as triangle inequality. Under this conditions, $(X,d)$ is said to be a metric
space.

Also, if $\Lambda$ is an abelian totally ordered group with a
total order $<$, a $\Lambda$-metric on a nonempty set $X$, is
defined as $d\colon X\times X\To\Lambda$ which satisfies
three above conditions, with this changing that 0, $<$ and $+$, are
interpreted to the identity element of $\Lambda$, the total
order on $\Lambda$, and the operation of group $\Lambda$,
respectively. Under this conditions, $(X,d)$ is named a
$\Lambda$-metric space. Since $(\mathbb{R},<)$, where $<$ is
ordinary order on $\mathbb{R}$, is an additive abelian totally
ordered group, each ordinary metric space is a $\Lambda$-metric
space. Of course, all notions and properties of the metric
spaces is not generalizable to $\Lambda$-metric spaces (see [2]).

In a $\Lambda$-metric space $X$, for every
$\varepsilon(>0)\in\Lambda$ and every $x\in X$, the open
$\varepsilon$-ball with center $x$, shown by
$B(x,\varepsilon)$, is defined as $\{y\in X\mid
d(x,y)<\varepsilon\}$, and the definition of limit
point, interior point, open and closed sets, and some another
relevent notions, are counterparts of the
definitions of same notions in ordinary metric spaces. The
collection of all open sets is formed the metric topology,
$\mathcal{T}_{d}$, on $X$, in each metric space $(X,d)$, and the
topological space $(X,\mathcal{T}_{d})$ is Hausdorff and normal.
The generalization of these, is established for any
$\Lambda$-metric space $(X,d)$ (See [2]).

We replace $P(S)$ (the power set of S), with $
\Lambda$ in definition of $\Lambda$-metric, where $S$ is a nonempty
set, of course with consideration a nonempty subset $Z$ of
$S$, and define a notion named $P_{Z}(S)$-metric and from that,
$P_{Z}(S)$-metric space, and then we make counterparts of some existent
notions in $\Lambda$-metric spaces for $P_{Z}(S)$-metric spaces, and in continuance, we proceed to make
counterparts of results, for $P_{Z}(S)$-metric spaces. Specially, we prove that the set of all open sets in every $P_{Z}(S)$-metric space $X$ is a topology for $X$, and that all open $\varepsilon$-balls in $X$ form a bases for this topology. 

The $\Lambda$-metrics, have fined plenty of
applications in group theory, especially in action of a group
on a set (See [1], [2], [3], [5]). Therefore it is possible that
$P_{Z}(S)$-metrics will be helpful in group theory as well general topology.

\maketitle

\section{Basic definitions and examples}

In the sequel we use [4], for standard terminology and notation on topology.

\

Before all, We give definition of $P_{Z}(S)$-metric.

\subsection{Definition}Let $P(S)$ be the power set of a
nonempty set $S$ and $\emptyset\neq Z\subseteq S$, and let $X$ be a
nonempty set. $(d,Z)$ is said to be a $P_{Z}(S)$-metric on $X$, if
$d\colon X\times X\To{P(S)}$, be a function with three following
properties:

1) $\forall x,y\in X,\ d(x,y)=\emptyset\Longleftrightarrow x=y$,

2) $\forall x,y\in X,\ d(x,y)=d(y,x)$,

3) $\forall x,y,z\in X,\ d(x,y)\subseteq d(x,z)\cup d(z,y) $.

\

We name the second property, symmetry, and the third,
 triangle inequality.

\

If there is no chance for ambiguity, we mention merely $P_{Z}(S)$-metric
 d, instead of $P_{Z}(S)$-metric $(d,Z)$.

\

The role of the set $Z$ will be specified in continuance.

\subsection{Definition}We say $(X,(d,Z))$ to be a $P_{Z}(S)$-metric
space, if $(d,Z)$ be a $P_{Z}(S)$-metric on $X$.

\

Again, where there is no chance for ambiguity, we mention only $P_{Z}(S)$-metric
 space $X$ instead of $P_{Z}(S)$-metric space $(X,(d,Z))$.

\subsection{Examples}a) Let $S$ be a nonempty set and let $Z$, $L$ be nonempty subsets of S,
and let $X$ be an arbitrary nonempty set. If we define:
\[d\colon X\times X\To{P(S)}\]
\[\ \ \ \ \ \ \ \ \ \ \ \ d(x,y)=\left\lbrace
\begin{array}{c l}
\emptyset&\text{\ if x=y},\\
L&\text{\ otherwise},
\end{array}
\right. \] then $(d,Z)$ is a $P_{Z}(S)$-metric on $X$ and therefore
$(X,(d,Z))$ is a $P_{Z}(S)$-metric space.

b) Let $X=\mathbb{R}, S=\mathbb{R}, Z=(a,b)$, where
$a,b\in\mathbb{R}$. If
\[d\colon X\times X\To{P(S)}\]
\[\ \ \ \ \ \ \ \ \ \ \ \ \ \ \ \ \ d(x,y)=\left\lbrace
\begin{array}{c l}
\emptyset&\text{\ if x=y},\\
\{x,y\}&\text{\ otherwise,}
\end{array}
\right.\]
then $(d,Z)$ is a $P_{Z}(S)$-metric on $X$.

\subsection{Definition}Let $(X,d)$ be a $P_{Z}(S)$-metric
space. If $\varepsilon\in P(S)$ contains $Z$, we say that
$\varepsilon$
is positive under $Z$, and we write $\varepsilon>_{Z}0$ or
$0<_{Z}\varepsilon$. If there is no chance for ambiguity, we
are satisfied with writing $\varepsilon>0$ or $0<\varepsilon$, respectively.

It is clear that in each $P_{Z}(S)$-metric space, $\emptyset$ is never  positive, and $S$, $Z$ are constantly
positive. Also, if $\varepsilon>0$ and $\varepsilon\subseteq\beta$, then $\beta>0$, where $\beta\in P(S)$.

\subsection{Definition}Let $(X,d)$ be a $P_{Z}(S)$-metric space and let $\varepsilon(\subseteq S)>0$,
and $x\in X$.
The open $\varepsilon$-ball with center $x$, is:
\[\{y\in X\mid d(x,y)\subset\varepsilon\},\]
and is denoted by $B(x,\varepsilon)$.

\

Under above conditions, $B(x,\varepsilon)-\{x\}$ is said to be deleted open $\varepsilon$-ball with center x
and is denoted by $B^{\ast}(x,\varepsilon)$.

\subsection{Examples}In the Example (2.3) in part (a) for all $x\in X$,
\[B(x,\varepsilon)=\left\lbrace
\begin{array}{c l}
X&\text{\ if $L\subset\varepsilon$},\\
\{x\}&\text{\ otherwise,}
\end{array}
\right.\]
and therefore if $L\subset Z$, then for all
$\varepsilon>0$ and all $x\in X$, $B(x,\varepsilon)=X$.

In part (b), for all $x\in X$,
\[B(x,\varepsilon)=\left\lbrace
\begin{array}{c l}
\varepsilon&\text{\ if $x\in\varepsilon$},\\
\{x\}&\text{\ otherwise,}
\end{array}
\right.\]

\

As it is observed in the above examples, in each
$P_{Z}(S)$-metric space, $x\in B(x,\varepsilon)$, for all
$\varepsilon>0$, the same that its similar is established in
$\Lambda$-metric spaces.

\subsection{Definition}Let $X$ be a $P_{Z}(S)$-metric space,

1)$x\in X$ is said to be an interior point of a subset $A$ of $X$, if,
\[\exists\varepsilon(\subseteq S)>0,\ B(x,\varepsilon)\subseteq A.\]
We denote by $A^{0}$ or $Int(A)$, the set of all interior points of
$A$ and name it inter of $A$.

2)$x\in X$ is said to be a limit point of a subset $A$ of $X$, if,
\[\forall\varepsilon(\subseteq S)>0,\ B^{\ast}(x,\varepsilon)\cap
A\neq\emptyset.\]
We denote by $A^{'}$, the set of all limit points of $A$.

3) A subset $U$ of $X$, is said to be open, if each its point be an interior point.

4)A subset $A$ of $X$, is said to be closed, if each its limit point, be an element of it.

\

It is observed that these definitions, are counterparts of the
definitions of similar notions in $\Lambda$-metric spaces.

\subsection{Remark}Clearly in each $P_{Z}(S)$-metric space, $\emptyset$ is open, by negation of antecedent and $X$ is
open by definition.

\subsection{Examples for open sets}In part (a) of (2.3), upon (2.6) if $A\subset X$, then $x\in A$
will be an interior point of $A$, iff there exist some
$\varepsilon>0$ so that $L$ be not a real subset of $\varepsilon$
and also, $X$ is open clearly. Therefore if $L\subset Z$, then all
$\gamma(>0)$, satisfy $L\subset\gamma$, and therefore excluding
$X$, $\emptyset$, there exist no another open sets, and the interior of
each another set, is $\emptyset$. But if there exist some member of $L$,
as $l$, so that it is not belong to $Z$, or $Z\subseteq L$, then each
 subset $A$ of $X$, is open, because in the first case, if $\varepsilon=Z$, and in the
second case, if $\varepsilon =L$, then  L is not a real subset of
$\varepsilon$. Therefore in that example, each set is open, or the open sets are restricted to $\emptyset$ and $X$.

In part (b) of (2.3), for each subset $A$ of $X$, upon (2.6),
$x\in A$ is an interior point of $A$, iff $x\notin(a,b)$ or
$x\in(a,b)\subseteq A$ (In the first case, take
$\varepsilon=(a,b)$, and in the second, $\varepsilon=A$).
Therefore, the open sets are only $X$, $\emptyset$ and each set $A$ that
$A\cap(a,b)=\emptyset$ or $(a,b)\subseteq A$.

\maketitle

\section{Statement of results}

In each metric space, every open $\varepsilon$-ball is an open
set. The following Lemma, is the counterpart of this case, for
$P_{Z}(S)$-metric spaces.

\subsection{Lemma}In every $P_{Z}(S)$-metric space, each open
$\varepsilon$-ball is an open set.

Proof: Let $B(x,\varepsilon)=U$, be an open $\varepsilon$-ball in
$P_{Z}(S)$-metric space $X$, and let $y\in U$. Set:
\[\delta=(\varepsilon-d(x,y))\cup Z.\]
Then $\delta>0$ and $B(y,\delta)\subseteq U$, since if $a\in
B(y,\delta)$, then,
\[d(a,x)\subseteq d(a,y)\cup
d(y,x)\subset\delta\cup d(x,y)=\varepsilon.\ \ \triangle\]

\subsection{Theorem}The collection of all open sets in a $P_{Z}(S)$-metric space
$X$, is a topology on $X$.

Proof: Let:
\[\mathcal{T}=\{U\subseteq X\mid U\ is\ open\ in\ X.\}.\]
We have:

1)Upon what we said above, $\emptyset, X\in\mathcal{T}.$

2)If $\{U_{i}\}_{i\in I}$, be a collection of elements of
$\mathcal{T}$, then $\bigcup_{i\in I}U_{i}$ is an open set in X,
and therefore is an element of $\mathcal{T}$, because if $x\in
\bigcup_{i\in I}U_{i}$, then:
\[\exists i_{o}\in I, \ x\in U_{i_{0}},\]
and since $U_{i_{0}}$ is open, then:
\[\exists i_{o}\in I,\ \ \exists \varepsilon_{i_{o}}({\subseteq X})>0, \  B(x,\varepsilon_{i_{o}})\subseteq U_{i_{0}}\subseteq\bigcup_{i\in I}U_{i}.\]

3)Let $\{U_{i}\}_{i=1}^n $ is a finite collection of members of
$\mathcal{T}$. Then $\bigcap_{i=1}^n U_{i}$ is an open set in X
and therefore is an element of $\mathcal{T}$, because if
$x\in\bigcap_{i=1}^n U_{i}$, then since each $U_{i}$ is open,
\[\forall i\in I,\,\ \exists\varepsilon_{i},\ B(x,\varepsilon_{i})\subseteq U_{i}.\]
Now set:
\[\varepsilon=\bigcap_{i=1}^n\varepsilon_{i}.\]
Since for each i, $\varepsilon_{i}>o$, then $\varepsilon>o$. Now,
\[\forall\ 1\leq i\leq n,\ B(x,\varepsilon)\subseteq B(x,\varepsilon_{i})\subseteq U_{i}.\]
Then $B(x,\varepsilon)\subseteq\bigcap_{i=1}^n U_{i}.\ \
\triangle$

\

Similar to $\Lambda$-metric spaces, the topology obtained from open sets of $P_{Z}(S)$-metric space $(X,d)$, is said to be
the $P_{Z}(S)$-metric topology obtained from $P_{Z}(S)$-metric $d$
on $X$, and is denoted by $\mathcal{T}_{d}$.

\subsection{Example}In the Example (a) of (2.3), if $L$ is not a real subset of
$Z$, the $P_{Z}(S)$-metric topology obtained from
$P_{Z}(S)$-metric d, is the discrete topology and else, that is
the indiscrete topology.

\subsection{Theorem}The collection of all open balls in each
$P_{Z}(S)$-metric space $(X,d)$, forms a bases for the
$P_{Z}(S)$-metric topology $\mathcal{T}_{d}$ on $X$.

Proof: Set:
\[\mathcal{B}=\{B(x,\varepsilon)\mid \varepsilon >0, x\in X\}.\]
Firstly, according to (3.1), the members of $\mathcal{B}$ are
open. Secondly, if $U$ is an open set in
$\mathcal{T}_{d}$, then for each $x\in U$, we have:
\[\exists \varepsilon_{x}>0,\ B(x,\varepsilon_{x})\subseteq U,\]
and since for each $x\in U$, $x\in B(x,\varepsilon_{x})$, then,
\[U=\bigcup_{x\in U}B(x,\varepsilon_{x}),\]
and therefore $\mathcal{B}$ is a bases for the topology
$\mathcal{T}_{d}$ on X.$\ \ \triangle$

\

Upon above theorem, to be an element of $X$ as an interior or limit
point for a set, and therefore to be a set, closed in a
$P_{Z}(S)$-metric space $(X,d)$, with meaning which we said
in (2.7), is equivalent with the analogous cases in topological
space $(x,\mathcal{T}_{d})$.

\

With the upper descriptions, if $(X,d)$, be a $P_{Z}(S)$-metric
space, then $(X,\mathcal{T}_{d})$ is a topological space and from
this, all contexts and corollaries, concerned with topological
spaces, are current in it.

\

If $(X,d)$ is a $P_{Z}(S)$-metric space and $(x_{n})$ is a
sequence in $X$, and if $x\in X$, then upon (3.4), $\lim_{
x\rightarrow\infty} x_{n}=x$ (in topological mining), iff,
\[\forall \varepsilon>0,\ \ \ \exists N\in\mathbb{N},\ \ \forall n\geq N,\ d(x_{n},x)\subset\varepsilon.\]

\subsection{Definition}Let $(X,d)$ be a $P_{Z}(S)$-metric space. A
sequence $(x_{n})$ in $X$ is said to be a Cauchy sequence, if,
\[\forall \varepsilon>0,\ \ \ \exists N\in\mathbb{N},\ \ \forall i,j\geq N,\ d(x_{i},x_{j})\subset\varepsilon.\]

\subsection{Theorem}Let $(X,d)$ be a $P_{Z}(S)$-metric space and $(x_{n})$ be a convergent sequence in
$X$. Then $(x_{n})$ is a Cauchy sequence in $X$.

Proof: Suppose that $\varepsilon>0$ and let $x_{n}\rightarrow x$,
where $x\in X$. We have,
\[\exists N\in\mathbb{N},\ \ \forall n\geq N,\ d(x_{n},x)\subset\varepsilon.\]
Now, if $i,j\geq N$, then,
\[d(x_{i},x)\subset\varepsilon,\ \ \ \ \ \ \ \  d(x_{j},x)\subset\varepsilon.\]
then,
\[d(x_{i},x_{j})\subseteq d(x_{i},x)\cup d(x,x_{j})\subset\varepsilon\cup
\varepsilon=\varepsilon.\ \ \triangle\]

\subsection{Theorem}If a Cauchy sequence $(x_{n})$ in $P_{Z}(S)$-metric space $(X,d)$,
 has a subsequence convergent to $x$, then $(x_{n})$ is convergent to $x$.

Proof: Let $(x_{n_{k}})$ be a subsequence of $(x_{n})$ such that
$x_{n_{k}}\rightarrow x$. With hypothesis for $\varepsilon>0$, we
have,
\[\exists N_{1}\in\mathbb{N},\ \ \forall k\geq N_{1},\ d(x_{n_{k}},x)\subset\varepsilon,\]
\[\exists N_{2}\in\mathbb{N},\ \ \forall i,j\geq N_{2},\ d(x_{i},x_{j})\subset\varepsilon.\]
Now if we take $N$ such that $N\geq N_{1}$ and $n_{N}\geq N_{2}$,
we have,
\[d(x_{i},x)\subseteq d(x_{i},x_{n_{i}})\cup d(x_{n_{i}},x)\subset\varepsilon\cup
\varepsilon=\varepsilon,\]
 for every $i\geq N$.$\ \ \triangle$

\subsection{Theorem}Let $(X_{1},d_{1})$ is a $P_{Z_{1}}(S_{1})$-metric
space and $(X_{2},d_{2})$ is a $P_{Z_{2}}(S_{2})$-metric
space. Then $f:X_{1}\rightarrow X_{2}$ is continuous (in
topological mining), iff for each $a\in X_{1}$, we have,
\[\forall \varepsilon>_{Z_{2}}0,\ \ \ \exists \delta>_{Z_{1}}0,\ \ \forall x\in X_{1},\ d_{1}(x,a)\subset\delta\Rightarrow d_{2}(f(x),f(a))\subset\varepsilon. \ \ (\ast)\]

Proof: Let $f:X_{1}\rightarrow X_{2}$ is continuous and let $a\in
X_{1}$ and $\varepsilon>_{Z_{2}}0$. Then upon definition of
the continuity, $f^{-1}(B(f(a),\varepsilon))$ is open in X.
Therefore since $a \in f^{-1}(B(f(a),\varepsilon))$,
\[\exists \delta>_{Z_{1}}0,\ B(a,\delta)\subseteq f^{-1}(B(f(a),\varepsilon)).\]
Now, if $x\in X_{1}$ and $d_{1}(x,a)\subset\delta$, then $x \in
f^{-1}(B(f(a),\varepsilon))$ and therefore
\[d_{2}(f(x),f(a))\subset\varepsilon.\]
Conversely, let $(\ast)$ is established. We prove that f is
continuous. Upon (3.4), it is enough that we prove that the
inverse image of each open $\varepsilon$-ball in $X_{2}$, is open
in $X_{1}$. Let $B(y,\varepsilon)$ be an open $\varepsilon$-ball
in $X_{2}$, and let $a \in f^{-1}(B(y,\varepsilon))$. Then upon
$(\ast)$,
\[\exists \delta>_{Z_{2}}0,\ \ \forall x\in X_{1},\ d_{1}(x,a)\subset\delta\Rightarrow d_{2}(f(x),f(a))\subset\varepsilon.\ \ (I)\]
We have $B(a,\delta)\subseteq f^{-1}(B(y,\varepsilon))$, because
if $x\in B(a,\delta)$, then $d_{1}(x,a)\subset\delta$ and
therefore upon (I), $d_{2}(f(x),f(a))\subset\varepsilon$. Then,
\[d_{2}(f(x),y)\subseteq d_{2}(f(x),f(a))\cup d_{2}(f(a),y)\subset\varepsilon\cup\varepsilon=\varepsilon.\ \ \triangle\]
\

At the end, we notice that for a nonempty set S, $P(S)$ with the
operation $\cup$, is not a group. Also, if $S(\neq\emptyset)$ is not
of one element, $(P(S),\subseteq)$ is not a totally ordered set.

We hope this paper give a motivation for characterization compact and Hausdorf $P_{Z}(S)$-metric spaces as well to prove a version of fixed point theorem for these spaces.

\maketitle

\[\mathbf{References}\]

[1] N. Brady, L. Ciobanu, A. Martino, S. O'Rourke \textit{The equation $x^{p}y^{q}=z^{r}$ and groups that act freely on $\Lambda$-trees,} Trans. Amer. Math. Soc. \textbf{361}, (2009) no. 1, 223-236.

[2] I. M. Chiswell, \textit{Introduction to $\Lambda$-trees,} World
Scientific Publishing Co. Inc. River Edge, NJ, 2001. 

[3] O. Kharlampovich, A. Myasnikov, and D. Serbin, \textit{Group action freely on $\Lambda$-trees,} arXiv: 0911.0209v4.
277-280

[4] J. R. Munkres, \textit{Topology: A First Course,} Prentice-Hall, Inc., Englewood Cliffs, N. J., 1975.

[5] S. O Rourke, \textit{Affine action on non-archimedean trees,} arXiv: 1112.4832v2 [math. GR].

\end{document}